\documentclass[12pt]{article}
\usepackage{amsfonts, amssymb, amsmath}
\usepackage{enumerate}
\usepackage{esint}

\newtheorem{thm}{Theorem}[section]

\newtheorem{lem}[thm]{Lemma}

\newtheorem{defn}[thm]{Definition}
\newtheorem{rem}[thm]{Remark}

\newcommand{\f}{\frac}

\newcommand{\vc}{\infty}

\newcommand{\su}{\subset}

\begin{document}

\title{ Weighted norm inequalities for spectral multipliers without Gaussian estimates}
\author{The Anh Bui\thanks{The Anh Bui was supported by a Macquarie University scholarship  \newline
{\it {\rm 2010} Mathematics Subject Classification:} 42B20, 42B35.
}}

\date{}


\maketitle
\begin{abstract}
Let $L$ be a non-negative self-adjoint operator on
$L^2(\mathbb{R}^n)$. By spectral theory, we can define the operator
$F(L)$, which is bounded on $L^2(X)$, for any bounded Borel function $F$. In this paper, we study the sharp weighted $L^p$ estimates for spectral multipliers
$F(L)$ and their commutators $[b, F(L)]$ with BMO functions $b$. We would like to emphasize that the Gaussian upper bound
condition on the heat kernels associated to  the semigroups
$e^{-tL}$ is not assumed in this paper.
 \end{abstract}
 \tableofcontents
\section{Introduction}
Suppose that
  $L$ is  a  non-negative self-adjoint
     operator   on $L^2({\mathbb{R}^n})$. Let $E(\lambda)$  be the spectral resolution
      of  $L$. By the spectral theorem, for any bounded Borel function \linebreak $m: [0, \infty)\rightarrow {\Bbb C}$,
one can define the operator
$$
F(L)=\int_0^{\infty} F(\lambda) dE(\lambda),
$$
 \noindent
 which is  bounded on $L^2(\mathbb{R}^n)$.\\

 The problem concerning the boundedness of $F(L)$ has attracted a lot of attention
and has been studied many authors, see for example \cite{A1, A2, B, C, MM, H, He, DeM, DOS}  and the references
therein.\\
\medskip

In many previous papers, the Gaussian upper bound condition plays an essential role, see for example \cite{A1, DeM, DOS}. Recall that the semigroup $\{e^{-tL}\}_{t>0}$ generated by
$L$ has the kernels $p_t(x,y)$ satisfy the Gaussian upper bounds if there exist $c, C>0$ so that
\begin{equation}\label{Gaussiansetimates}
|p_t(x,y)|\leq \f{C}{t^{n/2}}\exp\Big(-c\f{|x-y|^2}{t}\Big)
\end{equation}
for all $t>0$ and $x,y\in \mathbb{R}^n$. It was proved in \cite{DOS} that if the bounded
Borel function $F:[0,\vc)\rightarrow \mathbb{C}$ satisfies the
following for some $s>n/2$
$$
\sup_{t>0}\|\eta\delta_tF\|_{W^\vc_s}<\vc
$$
where $ \delta_t F(\lambda)=F(t\lambda)$,
   $\| F \|_{W^p_s}=\|(I-d^2/d x^2)^{s/2}F\|_{L^p}$ and $\eta$ is an auxiliary non-zero cut-off
    function such that
    $\eta \in C_c^{\infty}(\mathbb{R}_+)$,
 then the spectral multiplier $F(L)$ is of weak type
$(1,1)$ and hence by duality arguments, $F(L)$ is bounded on
$L^p$ for all $p\in (1,\vc)$. Then, The weighted $L^p$ estimates for spectral
multipliers $F(L)$ have been studied in \cite{DSY, A} recently.
However, there are many important operators $L$ which do not satisfy
(\ref{Gaussiansetimates}). It is natural to raise a question of the
boundedness of spectral multipliers $F(L)$ without the Gaussian upper bound condition (\ref{Gaussiansetimates}). The positive answers were addressed in \cite{B, DP, DY}. In
\cite{B}, S. Blunck replaced the pointwise kernel bounds
(\ref{Gaussiansetimates}) by the $L^p-L^q$
estimates to obtain the boundedness of $F(L)$ on $L^p$ spaces for an
appropriate range of $p$ whenever the bounded Borel function $F$ satisfies $\sup_{t>0}\|\eta\delta_tF\|_{W^2_s}<\vc$ with $s>n/2+1/2$. Under the similar condition $\sup_{t>0}\|\eta\delta_tF\|_{W^2_s}<\vc$ with $s>n/2+1/2$, the authors in \cite{DP} shown that if the semigroup $\{e^{-tL}\}$ Davies-Gaffney estimates, then $F(L)$ bounded on $H^1_L$, the Hardy space associated to operator $L$. The slightly better results were obtained by Duong and Yan in \cite{DY} which we sketch out briefly here. Assume that a nonnegative self-adjoint operator $L$ generates the semigroup $\{e^{-tL}\}$ which satisfies Davies-Gaffney estimates. If the following condition holds for some $s>n/p-n/2$ with $0<p\leq 1$,
\begin{equation}\label{DY'scondition}
     \sup_{t>0}\| \eta \delta_t F \|_{C^s} < \infty,
\end{equation}
where
$$
\|m\|_{C^s}=\begin{cases} \sum_{k=0}^s \sup_{\lambda\in
\mathbb{R}}|m^{(k)}(\lambda)| \ &\text{if $s\in \mathbb{Z}$}\\
\sum_{k=0}^{[s]} \sup_{\lambda\in \mathbb{R}}|m^{(k)}(\lambda)| +
\|m^{([s])}\|_{{\rm Lip}(s-[s])}\ &\text{if $s\in \mathbb{Z}$}
\end{cases}
$$
 with $[s]$ is a integer part of $s$, then $F(L)$ is bounded on $H^p_L(X)$, the Hardy space
associated to the operator $L$.  Moreover, it was also proved that if $L$
generates the semigroup $e^{-tL}$ satisfying the $L^{q_0}-L^2$
off-diagonal estimates for some $q_0\in [1, 2)$ then $F(L)$ is bounded
on $L^p(w)$ for $2<p<q_0'$ and $w\in A_{p/2}\cap RH_{(q_0'/p)'}$ and hence by duality, $F(L)$ is bounded
on $L^p(w)$  for $q_0<p<2$ and $w\in A_{p/q_0}\cap RH_{(2/p)'}$.\\
\medskip

In this paper, we establish the sharp weighted estimates for
spectral multipliers $F(L)$ and the commutators $[b,F(L)]$ of $F(L)$
with BMO functions $b$, where $L$ generates the semigroup $e^{-tL}$
satisfying the $L^{q_0}-L^2$ off-diagonal estimates (see Section 2 for precise definition). Precisely, it has
been shown that if $L$ generates the semigroup $e^{-tL}$ satisfying
the $L^{q_0}-L^2$ off-diagonal estimates and the bounded Borel function $F$ satisfies the following condition for $s>n/2$
\begin{equation}\label{cond1}
\sup_{t>0}\|\eta\delta_tF\|_{W^\vc_s}<\vc,
\end{equation}
then there exists $q_0<r_0<2$ so that $F(L)$ and the
commutators $[b,F(L)]$ are bounded on $L^p(w)$ for $r_0<p<q_0'$ and
$w\in A_{p/r_0}\cap RH_{(q_0'/p)'}$ and hence by duality $F(L)$ and the
commutators $[b,F(L)]$ are bounded on $L^p(w)$ for $q_0<p\leq r_0$ and $w\in
A_{p/q_0}\cap RH_{(r_0'/p)'}$. Note that the class $A_{p/r_0}\cap RH_{(q_0'/p)'}$ is larger than  the class $A_{p/2}\cap RH_{(q_0/p)'}$ for $p>2$ and also larger than the class $A_{p/q_0}\cap RH_{(2/p)'}$ for $r_0<p<2$, and the class  $A_{p/q_0}\cap RH_{(r_0'/p)'}$ is larger than the class $A_{p/q_0}\cap RH_{(2/p)'}$ for $r_0<p<2$ for $q_0<p\leq r_0$. Hence,  the obtained results in our paper are better than those in those in \cite{DY}. Moreover, since $r_0<2$, we also obtain the weighted estimate of $F(L)$ and the
commutators $[b,F(L)]$ on weighted $L^p$ spaces when $p=2$. It seems that the obtained results in \cite{DY} did not tell us the weighted estimate of $F(L)$ on weighted $L^2$ spaces.\\
\medskip

The outline of this paper is as follows. In section 2, some basic
properties of Mukenhoupt's weights are recalled and then a criterion
on weighted estimates for singular integrals in \cite{AM1} is
addressed. Section 3 is dedicated to study the weighted estimated
for spectral multipliers $F(L)$ and the commutators $[b,F(L)]$.

\section{Mukenhoupt's weights and  weighted estimates for singular integrals}
\subsection{Mukenhoupt's weights}

Throughout this article, we will often just use $B$ for $B(x_B,
r_B):=\{x:|x-x_B|\leq r_B\}$. Also given $\lambda > 0$, we will
write $\lambda B$ for the $\lambda$-dilated ball, which is the ball
with the same center as $B$ and with radius $r_{\lambda B} = \lambda
r_B$. For each ball $B\subset \mathbb{R}^n$ we set
$$
S_0(B)=B \ \text{and} \ S_j(B) = 2^jB\backslash 2^{j-1}B \
\text{for} \ j\in \mathbb{N}.
$$

We shall denote $w(E) :=\int_E w(x)dx$ for any measurable set $E
\subset \mathbb{R}^n$. For $1 \leq p \leq \infty$ let $p'$ be the
conjugate exponent of $p$, i.e. $1/p + 1/p' = 1$.

We first introduce some notation. We use the notation
$$
\fint_B h(x)dx=\f{1}{|B|}\int_Bh(x)dx.
$$
A weight $w$ is a non-negative locally integrable function on $\mathbb{R}^n$.
We say that $w \in A_p$, $1 < p < \infty$, if there exists a
constant $C$ such that for every ball $B \subset \mathbb{R}^n$,
$$
\Big(\fint_B w(x)dx\Big)\Big(\fint_B w^{-1/(p-1)}(x)dx\Big)^{p-1}\leq C.
$$
For $p = 1$, we say that $w \in A_1$ if there is a constant $C$ such
that for every ball $B \subset \mathbb{R}^n$,
$$
\fint_B w(y)dy \leq Cw(x) \ \text{for a.e. $x\in B$}.
$$
The reverse H\"older classes are defined in the following way: $w
\in RH_q, 1 < q < \infty$, if there is a constant $C$ such that for
any ball $B \subset \mathbb{R}^n$,
$$
\Big(\fint_B w^q(x) dx\Big)^{1/q} \leq C \fint_B w(x)dx.
$$
The endpoint $q = \infty$ is given by the condition: $w \in
RH_\infty$ whenever, there is a constant $C$ such that for any ball
$B \subset \mathbb{R}^n$,
$$
w(x)\leq C \fint_B w(y)dy  \ \text{for a.e. $x\in B$}.
$$
Let $w \in A_p$, for $1\leq p <\infty$, the weighted spaces $L^p_w$
can be defined by
$$\Big\{f :\int_{\mathbb{R}^n} f(x)^p w(x)dx < \infty\Big\}$$
with the norm $$\|f\|_{L^p(w)}=\Big(\int_{\mathbb{R}^n} f(x)^p w(x)dx\Big)^{1/p}.$$\\

We sum up some of the properties of $A_p$ classes in the following
results, see \cite{Du}.
\begin{lem}\label{weightedlemma1}
The following properties hold:
\begin{enumerate}[(i)]
\item $A_1\subset A_p\subset A_q$ for $1\leq p\leq q\leq \infty$.
\item $RH_\infty \su RH_q \su RH_p$ for $1\leq p\leq q\leq \infty$.
\item If $w \in A_p, 1 < p < \vc$, then there exists $1 < q < p$ such that $w \in A_q$.
\item If $w \in RH_q, 1 < q < 1$, then there exists $q < p < 1$ such that $w \in RH_p$.
\item $A_\vc =\cup_{1\leq p<\vc}A_p = \cup_{1< p\leq \vc}RH_p$
\end{enumerate}
\end{lem}

\subsection{Weighted norm inequalities for singular integrals}

 \begin{thm}\label{thm-Martell}

Let $1< p_0< q_0\leq \infty.$  Let $T$ be a bounded sublinear
operator on $L^{p_0}(\mathbb{R}^n)$, Let $\{A_r\}_{r>0}$ a family of
operators acting on $L^{p_0}(\mathbb{R}^n)$. Assume that

\begin{eqnarray}\label{e1-Martell}
\Big( \fint_{B} \big| T(I-A_{r_B})f\big|^{p_0}dx\Big)^{1/p_0} \leq
C \sum_j\alpha_j \Big(\fint_{2^{j}B}|f|^{p_0} \Big)^{1/p_0}
\end{eqnarray}
and
\begin{eqnarray}\label{e2-Martell}
 \Big( \fint_{B} \big| TA_{r_B}f\big|^{q_0}dx\Big)^{1/q_0} \leq
C \sum_j\alpha_j \Big(\fint_{2^{j}B}|Tf|^{p_0} \Big)^{1/p_0}
\end{eqnarray}

\noindent for all $f \in L^{\infty}_c(\mathbb{R}^n) $, and all ball
$B$ with radius $r_B$.

  If $\sum_j j\alpha_j<\infty$, then for all
 $p_0<p<q_0$  and $w\in A_{p/p_0}\cap RH_{(q_0/p)'}$,
there exists a constant $C$ such that
\begin{eqnarray}\label{e3}
   \|T f\|_{L^p(\mathbb{R}^n, w)}\leq    C\|  f\|_{L^p(\mathbb{R}^n, w)}.
\end{eqnarray}
and
\begin{eqnarray}\label{e4}
   \|[b,T] f\|_{L^p(\mathbb{R}^n, w)}\leq    C\|b\|_{{\rm BMO}}\|  f\|_{L^p(\mathbb{R}^n, w)}
\end{eqnarray}
for all $b\in {\rm BMO}$.
\end{thm}
\emph{Proof:} The proof of this theorem is just a combination the
arguments of Theorems 3.7 and 3.16 in \cite{AM1} and we omit details
here.

\section{Weighted estimates for spectral multipliers}
Let $T$ be a bounded linear operator on $L^2(\mathbb{R}^n)$. By the
kernel $K_T(x,y)$ associated to $T$ we mean that
$$
Tf(x)=\int_{\mathbb{R}^n}K_T(x,y)f(y)dx
$$
holds for all $f\in L^\vc(\mathbb{R}^n)$ with compact support and
for all $x\notin $ supp$f$.

\begin{defn}[\cite{AM2}]
Let $1\leq p\leq q\leq \vc$. We say that the family $\{T_t\}_{t>0}$
of sublinear operators satisfies $L^p-L^q$ full off-diagonal
estimates,  in short $T_t\in \mathcal{F}(L^{p}-L^q)$, if there
exists some $c>0$, for all closed sets $E$ and $F$, all $f$ with
supp$f\subset E$ and all $t>0$ so that
\begin{equation}\label{eq-full-off-diag-estimates}
\|T_t f\|_{L^q(F)}\leq c
t^{-\f{n}{2}(\f{1}{p}-\f{1}{q})}\exp\Big(-c\f{d^2(E,F)}{t}\Big)\|f\|_{L^{p}}.
\end{equation}
\end{defn}

Let us summarize some basic properties concerning the classes
$\mathcal{F}(L^{p}-L^q)$, see \cite{AM2}.
\begin{enumerate}[(i)]

\item For $p\leq p_1\leq q_1\leq q$, $\mathcal{F}(L^{p_1}-L^{q_1})\subset
\mathcal{F}(L^{p}-L^{q})$.

\item $T_t\in \mathcal{F}(L^{1}-L^\vc)$ if and only if the
associated kernel $p_t(x,y)$ of $T_t$ satisfies the Gaussian upper
bound, that is, there exist positive constants $c$ and $C$ so that
$$
|p_t(x,y)|\leq \f{C}{t^{n/2}}\exp\Big(-c\f{|x-y|^2}{t}\Big)
$$
for all $x, y\in \mathbb{R}^n$ and $t>0$.

\item $T_t\in \mathcal{F}(L^{p}-L^q)$ if and only if
$T^*_t\in \mathcal{F}(L^{q'}-L^{p'})$.

\item If $S_t\in \mathcal{F}(L^{p}-L^r)$ and $T_t\in
\mathcal{F}(L^{r}-L^q)$ for $p\leq r\leq q$, then $T_t\circ S_t \in
\mathcal{F}(L^{p}-L^q)$.
\end{enumerate}

Full off-diagonal estimates appear when dealing with semigroups of
second order elliptic operators (see for example \cite{LSV, Aus}) or
semigroups of Shr\"odinger operators with real potentials \cite{As}.
The most studied case is when $p = 1$ and $q = \vc$ which means that
the kernel of $T_t$ has pointwise Gaussian upper bounds.

Let $L$ be a non-negative self-adjoint operator. Assume that the
semigroup $e^{-tL}$, generated by $-L$ on $L^2(X)$ satisfies
$L^{q_0}-L^2$ full off-diagonal estimates for some $q_0\in [1,2)$.
By (iii), $e^{-tL}\in L^{2}-L^{q_0'}$. Since
$e^{-tL}=e^{-\f{t}{2}L}\circ e^{-\f{t}{2}L}$, by (iv) we have
$e^{-tL}\in \mathcal{F}(L^{q_0}-L^{q_0'})$. In this paper, we will work with operators which generate the semigroup $e^{-tL}\in L^{q_0}-L^{q_0'}$ for some $q_0\in [1,2)$.

The following result is the main result of this paper.
\begin{thm}\label{mainresult}
Let $L$ be a self-adjoint non-negative operator on $L^2$ satisfying
the $L^{q_0}-L^2$ full off-diagonal for some $q_0\in [1,2)$. Set  $r_0=\max(q_0, \f{n}{s})$. If the bounded
Borel function $F:[0,\vc)\rightarrow \mathbb{C}$ satisfies the
following for some $s>n/2$
$$
\sup_{t>0}\|\eta\delta_tF\|_{W^\vc_s}<\vc
$$
where $ \delta_t F(\lambda)=F(t\lambda)$,
   $\| F \|_{W^p_s}=\|(I-d^2/d x^2)^{s/2}F\|_{L^p}$ and $\eta$ is an auxiliary non-zero cut-off
    function such that
    $\eta \in C_c^{\infty}(\mathbb{R}_+)$,
 then

(a) $F(L)$ is bounded on $L^p(w)$ for all $r_0<p<q_0'$ and $w\in
A_{p/r_0}\cap RH_{(q_0'/r_0)'}$;

(b) moreover, for $b\in BMO(\mathbb{R}^n)$, the commutator $[b,
F(L)]$ is also bounded on $L^p(w)$ for all $r_0<p<q_0'$ and $w\in
A_{p/r_0}\cap RH_{(q_0'/r_0)'}$.
\end{thm}

Before coming to the proof of Theorem \ref{mainresult} we would like to compare the obtained results in  Theorem \ref{mainresult} with those in Theorem 5.2 in \cite{DY}. Under the assumption as in Theorem \ref{mainresult}, Theorem 5.2 in \cite{DY} implies that $F(L)$ is bounded on $L^p(w)$ for either $2<p<q_0'$ and $w\in A_{p/2}\cap RH_{(q_0'/p)'}$  or $q_0<p<2$ and $w\in A_{p/q_0}\cap RH_{(2/p)'}$. Since the class $A_{p/r_0}\cap RH_{(q_0'/p)'}$ is larger than  the class $A_{p/2}\cap RH_{(q_0/p)'}$ for $p>2$ and also larger than the class $A_{p/q_0}\cap RH_{(2/p)'}$ for $r_0<p<2$, and the class  $A_{p/q_0}\cap RH_{(r_0'/p)'}$ is larger than the class $A_{p/q_0}\cap RH_{(2/p)'}$ for $r_0<p<2$ for $q_0<p\leq r_0$,  the obtained results in Theorem \ref{mainresult} are better than those in \cite{DY}. Moreover, since $r_0<2$, we also obtain the weighted estimate of $F(L)$ and the
commutators $[b,F(L)]$ on weighted $L^p$ spaces when $p=2$. It seems that the obtained results in \cite{DY} did not tell us the weighted estimate of $F(L)$ on weighted $L^2$ spaces.\\

\medskip

We split the proof of Theorem \ref{mainresult} into a few lemmas. 

\begin{lem}\label{lem1}
Let $p\in (q_0,2)$ and $F$ be a bounded Borel function with {\rm supp} $F\subset [0, R]$. There exists a constant $C>0$ such that
$$
\|F(\sqrt{L})f\|_{L^2(B)}\leq CR^{n(1/p-1/2)}\|f\|_{L^p(S_j(B))}\|F\|_{L^\vc}
$$
for all balls $B$, $j\geq 3$ and all $f\in L^p(S_j(B))$.
\end{lem}
\emph{Proof:} Setting $G(\lambda)=e^{\lambda^2/R^2}F(\lambda)$, then $\|G\|_{L^\vc}\approx \|F\|_{L^\vc}$. Moreover, we have $F(\sqrt{L})f=G(L)e^{-\f{1}{R^2}L}f$. Therefore,
\begin{equation*}
\begin{aligned}
\|F(\sqrt{L})f\|_{L^2(B)}&=\|G(L)e^{-\f{1}{R^2}L}f\|_{L^2(B)}\leq \|G(L)e^{-\f{1}{R^2}L}\|_{L^p(S_j(B))\rightarrow L^2(B)} \|f\|_{L^p(S_j(B))}\\
&\leq \|G(L)\|_{L^2(\mathbb{R}^n)\rightarrow L^2(\mathbb{R}^n)}\times \|e^{-\f{1}{R^2}L}\|_{L^p(\mathbb{R}^n)\rightarrow L^2(\mathbb{R}^n)}\times \|f\|_{L^p(S_j(B))}\\
&\leq CR^{n(1/p-1/2)}\|G\|_{L^\vc}\|f\|_{L^p(S_j(B))}\approx CR^{n(1/p-1/2)}\|f\|_{L^p(S_j(B))}\|F\|_{L^\vc}.
\end{aligned}
\end{equation*}

\begin{lem}\label{lem2}
For any $p\in (q_0, 2)$, there exist two constants $C>0$ and $c>0$ so that for all closed sets $E$ and $F$, all $f$ with {\rm supp} $f \subset E$ and all $z\in \mathbb{C}_+=\{z\in \mathbb{C}: \Re z>0\}$, the following holds
$$
\|e^{-zL}f\|_{L^2(F)}\leq C(|z|\cos\theta)^{-\f{n}{2}(\f{1}{p}-\f{1}{2})}\exp\Big(-c\f{d^2(E,F)}{|z|}\cos\theta\Big)\|f\|_{L^p(E)}
$$
where $\theta=\arg z$.
\end{lem}

To prove Lemma \ref{lem2}, we need the following version of Phragmen-Lindel\"of Theorem, see for example \cite[Lemma 9]{Da}.
\begin{lem}\label{PLlemma}
Suppose that function $G$ is analytic in  $\{z\in \mathbb{C}: \Re z>0\}$ and that
$$
|G(|z|e^{i\theta})|\leq a_1(|z|\cos \theta)^{-\beta_1},
$$
$$
|G(|z|)|\leq a_1 |z|^{-\beta_1}\exp(-a_2|z|^{-\beta_2})
$$
for some $a_1, a_2>0, \beta_1\geq 0, \beta_2\in (0,1]$, all $|z|>0$ and all $\theta\in (-\pi/2,\pi/2)$. Then
$$
|G(|z|e^{i\theta})|\leq 2^{\beta_1}a_1(|z|\cos \theta)^{-\beta_1}\exp\Big(-\f{a_2\beta_2}{2}|z|^{-\beta_2}\cos \theta\Big)
$$
for all $|z|>0$ and all $\theta\in (-\pi/2,\pi/2)$.
\end{lem}
\emph{Proof of Lemma \ref{lem2}: } Let $\|g\|_{L^2}=1$ supported in $F$. We define the holomorphic function $G_f: \mathbb{C}_+\rightarrow \mathbb{C}$ by setting
$$
G_{f}(z)=\int e^{-zL}f(x)g(x)dx.
$$
For any $z\in \mathbb{C}_+$, we have
\begin{equation*}
\begin{aligned}
G(z)\leq \|e^{-zL}f\|_{L^2(F)} &=\|e^{-i\Im z L}\circ e^{-\Re z L}f\|_{L^2(F)}\\
&\leq C\|e^{-i\Im z L}\|_{L^2(\mathbb{R}^n)\rightarrow L^2(\mathbb{R}^n)}\|e^{-\Re z L}\|_{L^p(\mathbb{R}^n)\rightarrow L^2(\mathbb{R}^n)}\|f\|_{_{L^p(E)}}\\
&\leq C(|z|\cos\theta)^{-\f{n}{2}(\f{1}{p}-\f{1}{2})}\|f\|_{_{L^p(E)}}.
\end{aligned}
\end{equation*}
In particular when $\theta=0$, we have
$$
G_{f}(|z|)\leq C(|z|)^{-\f{n}{2}(\f{1}{p}-\f{1}{2})}\exp\Big(-c\f{d^2(E,F)}{|z|}\Big)\|f\|_{_{L^2(E)}}.
$$
At this stage, applying Lemma \ref{PLlemma} with $a_1=C \|f\|_{L^p(E)}$, $a_2=-cd^2(E,F)$, $\beta_1=\f{n}{2}(1/p-1/2)$ and $\beta_2=1$, we obtain the desired estimate.

\begin{lem}\label{lem3}
Let $p\in (q_0,2)$ and $R>0, s>0$. For any $\epsilon >0$, there exists a constant $C=C(\epsilon,s)>0$ so that
\begin{equation}\label{eq1-lem3}
\|F(\sqrt{L})f\|_{L^2(B)}\leq C\f{R^{n(1/p-1/2)}}{(2^jr_B R)^s}\|\delta_R F\|_{W^\vc_{s+\epsilon}}\|f\|_{L^p(S_j(B))}
\end{equation}
for all balls $B$, all $j\geq 3$, all $f\in L^p(S_j(B))$ and all bounded Borel function $F$ supported in $[R/4, R]$.
\end{lem}
\emph{Proof:} Using the Fourier inversion formula, we write
$$
G(L/R^2)e^{-\f{1}{R^2}L}=c\int_{\mathbb{R}}e^{-\f{1-i\tau}{R^2}L}\widehat{G}(\tau)d\tau.
$$
Hence,
$$
F(\sqrt{L})f=c\int_{\mathbb{R}}\widehat{G}(\tau)e^{-\f{1-i\tau}{R^2}L}fd\tau
$$
where $G(\lambda)=[\delta_R F](\sqrt{\lambda})e^{\lambda}$.

Applying Lemma \ref{lem2}, we have, for any $f$ supported in $S_j(B)$,
\begin{equation*}
\begin{aligned}
\|F(\sqrt{L})f\|_{L^2(B)}&\leq c\int_{\mathbb{R}}\widehat{G}(\tau)\|e^{-\f{1-i\tau}{R^2}L}f\|_{L^2(B)}d\tau\\
&\leq cR^{n(1/p-1/2)}\int_{\mathbb{R}}\widehat{G}(\tau)\exp\Big(-c\f{(2^jr_BR)^2}{(1+\tau^2)}\Big) d\tau\times \|f\|_{L^p(S_j(B))}\\
&\leq cR^{n(1/p-1/2)}\|f\|_{L^p(S_j(B))}\int_{\mathbb{R}}\widehat{G}(\tau)\f{(1+\tau^2)^{s/2}}{(2^jr_BR)^s} d\tau\\
&\leq c\f{R^{n(1/p-1/2)}}{(2^jr_BR)^s}\|f\|_{L^p(S_j(B))}\Big(\int_{\mathbb{R}}|\widehat{G}(\tau)|^2(1+\tau^2)^{s+\epsilon+1/2} d\tau\Big)^{1/2}\\
&~~~~~\times \Big(\int_{\mathbb{R}}(1+\tau^2)^{-\epsilon-1/2} d\tau\Big)^{1/2}\\
&\leq c\f{R^{n(1/p-1/2)}}{(2^jr_BR)^s}\|G\|_{W^2_{s+\epsilon +1/2}} \|f\|_{L^p(S_j(B))}.
\end{aligned}
\end{equation*}
Note that since supp $F\subset [R/4, R]$, we have $\|G\|_{W^2_{s+\epsilon +1/2}}\leq C\|\delta_RF\|_{W^2_{s+\epsilon +1/2}}\leq C\|\delta_RF\|_{W^\vc_{s+\epsilon +1/2}}$, and so
\begin{equation}\label{eq2-lem3}
\|F(\sqrt{L})f\|_{L^2(B)}\leq c\f{R^{n(1/p-1/2)}}{(2^jr_BR)^s}\|\delta_RF\|_{W^\vc_{s+\epsilon +1/2}} \|f\|_{L^p(S_j(B))}.
\end{equation}
To get rid of $1/2$ on the RHS of (\ref{eq2-lem3}), we use the interpolation arguments as in \cite{MM, DOS}. We first note that (\ref{eq1-lem3}) is equivalent to the following estimate
\begin{equation*}
\begin{aligned}
\|\delta_{1/R}H(\sqrt{L})f\|_{L^2(B)}\times (2^jr_BR)^s\leq C R^{n(1/p-1/2)}\|f\|_{L^p(S_j(B))}\|H\|_{W^\vc_{s+\epsilon}}.
\end{aligned}
\end{equation*}
for all bounded Borel function $H$ with supp $H\subset [1/4,1]$.

Now we define the linear operator $\mathcal{A}_{R,f}: L^\vc([1/4,1])\rightarrow L^2(B,dx)$ by setting
$$
\mathcal{A}_{R,f}(H)=\delta_{1/R}H(\sqrt{L})f.
$$
By Lemma \ref{lem1},
$$
\|\mathcal{A}_{R,f}\|_{L^\vc([1/4,1])\rightarrow L^2(B,dx)}\leq CR^{n(1/2-1/p)}\|f\|_{L^p(S_j(B))}.
$$

Setting $d\mu_{s,R}=(R2^jr_B)^{s}dx$, then (\ref{eq2-lem3}) tells us that
$$
\|\mathcal{A}_{R,f}\|_{W^\vc_{s+1/2+\epsilon}([1/4,1])\rightarrow L^2(B,\mu_{s,R})}\leq CR^{n(1/p-1/2)}\|f\|_{L^p(S_j(B))}.
$$
By interpolation, for each $\theta\in (0,1)$ there exists a constant $C$ such that
$$
\|\mathcal{A}_{R,f}(H)\|_{L^2(B,\mu_{s\theta ,R})}\leq CR^{n(1/p-1/2)}\|f\|_{L^p(S_j(B))}\|H\|_{[L^\vc, W^\vc_{s+1/2+\epsilon}]_{[\theta]}}.
$$
Therefore, for all $s>0, \epsilon' >0$ and $\theta\in (0,1)$,
$$
\|\mathcal{A}_{R,f}(H)\|_{L^2(B,\mu_{s\theta ,R})}\leq CR^{n(1/p-1/2)}\|f\|_{L^p(S_j(B))}\|H\|_{W^\vc_{s\theta+\epsilon'+\theta/2}}.
$$
By choosing $s'=s/\theta$ and taking $\theta$ small enough we obtain
$$
\|\mathcal{A}_{R,f}(H)\|_{L^2(B,\mu_{s'\theta ,R})}\leq CR^{n(1/p-1/2)}\|f\|_{L^p(S_j(B))}\|H\|_{W^\vc_{s'+\epsilon''}}.
$$
This completes our proof.

\medskip

We are now in position to prove Theorem \ref{mainresult}.

\emph{Proof of Theorem \ref{mainresult}:} Take $p_0\in (r_0,  2)$.
  Let $M\in{\mathbb N}$ such that $M>s/2$. We will show that
  (\ref{e1-Martell}) and (\ref{e2-Martell}) hold for $T=m(\sqrt{L})$ and $A_{r_B}=I-(I-e^{-r_B^2L})^M$.
 To verify (\ref{e1-Martell}), we will show that  for all balls $B$,
\begin{equation}\label{e1-mainthm}
\Big( \fint_{B} \big| F(L)(I-A_{r_B})f\big|^{p_0}dx\Big)^{1/p_0}
\leq C \sum_j\alpha_j \Big(\fint_{2^{j}B}|f|^{p_0} \Big)^{1/p_0}
\end{equation}

\noindent
  for all $f\in L^{\infty}_c(X)$, where $\alpha_j=2^{-j(s-\f{n}{p_0})}$.

Let us prove (\ref{e1-mainthm}). Since
$\sup_{t>0}\|\phi(\cdot)F(t\cdot)\|_{W^\vc_s} \approx
\sup_{t>0}\|\phi(\cdot)\widetilde{F}(t\cdot)\|_{W^\vc_s}$ where
$\widetilde{F}(\lambda)=F(\sqrt{\lambda})$, instead of proving
(\ref{e1-mainthm}) for $F(L)$, we will prove (\ref{e1-mainthm}) for
$F(\sqrt{L})$.

Let $\phi_\ell$ denote the function $\phi(2^{-\ell}\cdot)$. Following the standard arguments, see for example \cite{DOS}, we can write
\begin{eqnarray*}
F(\lambda)=
\sum_{\ell=-\infty}^{\infty}\phi(2^{-\ell}\lambda)F(\lambda) =
\sum_{\ell=-\infty}^{\infty}F^{\ell}(\lambda), \ \ \  \ \forall\,
\lambda >0.
\end{eqnarray*}

  For every   $\ell\in{\Bbb Z}$ and $r>0$, we set for $\lambda>0,$

\begin{eqnarray*}
F_{r, M}(\lambda)&=& F(\lambda)(1-e^{-(r\lambda)^2})^{M},\\[5pt]
F^{\ell}_{r, M}(\lambda)&=&
F^{\ell}(\lambda)(1-e^{-(r\lambda)^2})^{M}.
\end{eqnarray*}

\noindent Given a  ball $B\subset X$, we write $
 f=\sum\limits_{j=0}^{\infty} f_j  $ in which $ f_j=f\chi_{S_j (B)}$.
We may  write

\begin{eqnarray*}
 F(\sqrt{L})(1-e^{- r_B^2L})^{M}f &=& F_{r_B, M}(\sqrt{L})f \\
 &= &
 \sum_{j=0}^2 F_{r_B, M}(\sqrt{L})f_j+\sum_{j=3}^\infty \sum_{l=-\infty}^\infty F^l_{r_B, M}(\sqrt{L})f_j.  \nonumber
\end{eqnarray*}

Since $e^{-tL}\in \mathcal{F}(L^{q_0}-L^{q_0'})$, we have  that for
any $t>0$, $\|e^{-tL}f\|_{L^p}\leq C\| f\|_{L^p}$ for all $p\in
(q_0,q_0')$. This,
 in combination with  $L^{p}$-boundedness of the operator $F(\sqrt{L})$
 (see Theorem 5.2, \cite{DY}), gives  that
 for all balls $B\ni x$,

\begin{eqnarray}\label{e4.9}
 \Big(  \fint_{B}
\big|F_{r_B, M}(\sqrt{L})f_j \big|^{p_0}dx\Big)^{1/p_0}
&\leq& |B|^{-{1/p_0}} \big\|F_{r_B, M}(\sqrt{L})f_j\big\|_{L^{p_0}(X)}\nonumber\\
&\leq&  C|B|^{-1/p_0} \big\|  f_j\big\|_{L^{p_0}(X)}\nonumber\\
&\leq& C  \Big(\fint_{2^jB}|f|^{p_0}\Big)^{1/p_0}
\end{eqnarray}

\noindent
  for  $j=1,2.$

  Fix $j\geq 3.$ Let $p_1\geq 2$. Since $p_0<2$, using H\"older's inequality, we have
 \begin{equation}\label{e4.10}
 \begin{aligned}
 \Big(  \fint_{B}
\big|F^{\ell}_{r_B, M}(\sqrt{L})  f_j \big|^{p_0}dx\Big)^{1/p_0}
 &\leq |B|^{ -{1\over 2}}  \big\|F^{\ell}_{r_B, M}(\sqrt{L})  f_j\big\|_{L^{2}(B)}.
\end{aligned}
\end{equation}
Note that supp $F^{\ell}_{r_B, M}\subset [2^{\ell -2}, 2^\ell]$. So, using Lemma \ref{lem3} we obtain, for $s>s'>n/p_0$,
$$
\big\|F^{\ell}_{r_B, M}(\sqrt{L})  f_j\big\|_{L^{2}(B)}\leq 2^{\ell n(1/p_0-1/2)}(2^\ell r_B)^{-s'}2^{-js'}\|\delta_{2^\ell}F^{\ell}_{r_B, M}\|_{W^\vc_{s}}\|f_j\|_{L^{p_0}}.
$$
Let $k$ be an integer so that $k>s$. Then we have
\begin{equation*}
\begin{aligned}
\|\delta_{2^\ell}F^{\ell}_{r_B, M}\|_{W^\vc_{s}}&\leq C\|\phi\delta_{2^\ell}F\|_{W^\vc_{s}}\|(1-e^{(2^\ell r_B\cdot)^2})^M\|_{C^k[1/4,1]}\\
&\leq C\sup_{t>0}\|\phi\delta_{t}F\|_{W^\vc_{s}}\min\{1, (2^\ell r_B)^{2M}\}.
\end{aligned}
\end{equation*}
Therefore,
$$
\big\|F^{\ell}_{r_B, M}(\sqrt{L})  f_j\big\|_{L^{2}(B)}\leq C2^{-s'j}\min\{1,(2^{\ell}r_B)^{2M}\}(2^\ell
r_B)^{-s'}2^{\ell n(\f{1}{p_0}-\f{1}{2})}\sup_{t>0}\|\phi\delta_{t}F\|_{W^\vc_{s}}\|f_j\|_{L^{p_0}}
$$
for all $\ell\in \mathbb{Z}$.

This together with (\ref{e4.10}) gives
\begin{equation*}
\begin{aligned}
\Big(  \fint_{B} &\big|F^{\ell}_{r_B, M}(\sqrt{L})  f_j
\big|^{p_0}dx\Big)^{1/p_0}\\
  &\leq 2^{-j(s'-n/p_0)}|B|^{{1\over p_0}-{1\over 2}} \min\{1,(2^{\ell}r_B)^{2M}\}(2^\ell
r_B)^{-s'}2^{\ell n(\f{1}{p_0}-\f{1}{2})}\\
&~~~~\times\sup_{t>0}\|\phi\delta_{t}F\|_{W^\vc_{s}}\Big(\fint^{2^jB}|f|^{p_0}\Big)^{1/p_0}\\
&\leq c2^{-j(s'-n/p_0)}\min\{1,(2^{\ell}r_B)^{2M}\}(2^\ell
r_B)^{-(s'-n/p_0+n/2)}\sup_{t>0}\|\phi\delta_{t}F\|_{W^\vc_{s}}\Big(\fint^{2^jB}|f|^{p_0}\Big)^{1/p_0}.
\end{aligned}
\end{equation*}
This follows (\ref{e1-mainthm}). The proof is completed.

Since $e^{-tL}\in \mathcal{F}(L^{q_0}-L^{q_0'})$,
$A_{r_B}\mathcal{F}(L^{q_0}-L^{q_0'})$. This together with the fact
that $A_{r_B}$ and $F(\sqrt{L})$ are commutative gives
\begin{eqnarray}
 \Big( \fint_{B} \big| F(\sqrt{L})A_{r_B}f\big|^{q'_0}dx\Big)^{1/q'_0} \leq
C \sum_j\alpha_j \Big(\fint_{2^{j}B}|F(\sqrt{L}) f|^{p_0}
\Big)^{1/p_0}.
\end{eqnarray}
Therefore, Theorem \ref{thm-Martell} tells us that $F(\sqrt{L})$ and
the commutator $[b, F(\sqrt{L})]$ are bounded on $L^p(w)$ for all
$p_0<p<q_0'$ and $w\in A_{p/p_0}\cap RH_{(q_0'/p_0)'}$. Since
$A_{p/r_0}=\cup_{p_0>r_0} A_{p/p_0}$ and
$RH_{(q_0'/r_0)'}=\cup_{p_0>r_0} RH_{(q_0'/p_0)'}$, letting
$p_0\rightarrow r_0$ we obtain the desired results.

This competes our proof.\\

\begin{rem}
(i) Since $r_0<2$, $r_0'>2>r_0$. By duality, the weighted $L^p$
estimates for $F(\sqrt{L})$ and the commutator $[b, F(\sqrt{L})]$
can be obtained for $q_0<p\leq r_0$. More precisely, $F(\sqrt{L})$
and the commutator $[b, F(\sqrt{L})]$ are bounded on $L^p(w)$ for
all $q_0<p\leq r_0$ and $w\in A_{p/q_0}\cap RH_{(r_0'/q_0)'}$.

(ii) It can be believed that the approach in this paper can be extended to spaces of homogeneous type. This will appear in the forth-coming paper, see for example \cite{AD}.
\end{rem}

\emph{Acknowledgements} \  \ The author would like to thank his supervisor, Prof. X. T. Duong for helpful
comments and suggestions.

\medskip

\noindent 
The Anh Bui\\
Department of Mathematics, Macquarie University, NSW 2109, Australia and \\
Department of Mathematics, University of Pedagogy, Ho chi Minh city, Vietnam \\
Email: the.bui@mq.ed.au and bt\_anh80@yahoo.com

\begin{thebibliography}{99}
\bibitem[An1]{A} B. T. Anh, Sharp weighted $L^p$ estimates for spectral
multipliers, unpublished note.


\bibitem[AD]{AD} B. T. Anh and X. T. Duong, Weighted norm inequalities for some singular integrals associated to Schr\"odinger operators with real potential on manifolds, in preparation.

 \bibitem[A1]{A1} G. Alexopoulos, Spectral  multipliers on Lie groups of polynomial growth,
 {\it Proc. Amer. Math. Soc.}, {\bf 46 } (1994), 457-468.

\bibitem[A2]{A2}
G.~Alexopoulos, Spectral multipliers for Markov chains.
 {\it J. Math. Soc. Japan}, {\bf 56}(3)(2004), 833--852.

\bibitem[As]{As} J. Assaad, Riesz transforms associated to
Schr\"odinger operators with negative potential, preprint.

\bibitem[Aus]{Aus} P. Auscher, On necessary and sufficient conditions for $L^p$
estimates of Riesz transform associated elliptic operators on $\mathbb{R}^n$ and
related estimates, \emph{Mem. Amer. Math. Soc.} \textbf{186} (2007), no. 871.

 \bibitem[AM1]{AM1} P. Auscher, J.M. Martell, Weighted norm inequalities, off-diagonal estimates and elliptic
operators. Part I: General operator theory and weights,
  {\it Adv. in Math.}, {\bf 212}(2007), 225-276.

 \bibitem[AM2]{AM2} P. Auscher, J.M. Martell, Weighted norm inequalities, off-diagonal estimates and elliptic
operators. Part II: Off-diagonal estimates on spaces of homogeneous
type,  Math. Z. 260 (2008), 527-539.


\bibitem[B]{B} S. Blunck,   A H\"ormander-type spectral multiplier theorem for operators
without heat kernel, {\it Ann. Sc. Norm. Super. Pisa Cl. Sci.}, {\bf
2} (2003),  449-459.

\bibitem[C]{C} M. Christ, $L^p$ bounds for spectral multipliers on nilpotent groups,
 {\it Trans. Amer. Math. Soc.}, {\bf 328} (1991), 73-81.







\bibitem[CSo]{CSo} M. Christ and C.D. Sogge, The weak type $L^1$ convergence
of eigenfunction expansions for pseudodifferential operators,
 {\it Invent. Math.}, {\bf 94} (1988), 421-453.



\bibitem[Da]{Da} E. B. Davies, Uniformly elliptic operators with measureable coefficients, J. Funct. Anal. 132 (1995), 141-169.


\bibitem[DS]{DS} E.B. Davies and B. Simon, $L^p$ norms of
non-critical Shr\"odinger semigroups, J. Funct. Anal. 102, 95-115,
(1991).



\bibitem[D]{Du} J. Duoandikoetxea, Fourier Analysis, Grad. Stud. math, 29, American Math. Soc., Providence,
2000.



\bibitem[DeM]{DeM}   L. De Michele and G. Mauceri, $H^p$ multpliers on stratified groups,
  {\it Ann. Mat. Pura Appl.},
{\bf 148} (1987), 353--366.

\bibitem[DOS]{DOS}  X.T. Duong, E.M. Ouhabaz  and A. Sikora,
 Plancherel-type estimates and sharp spectral multipliers.
{\it J. Funct. Anal.}, {\bf 196} (2002),  443-485.

\bibitem[DY]{DY} X.T. Duong, and L.X. Yan, Spectral multipliers for
Hardy spaces associated to non-negative self-adjoint operators
satisfying Davies-Gaffney estimates, to appear in J. Math. Soc.
Japan.

\bibitem[DSY]{DSY} X.T. Duong, A. Sikora and L. Yan,
Weighted norm inequalities, Gaussian bounds and sharp spectral
multipliers, to appear J. Funct. Anal..

\bibitem[DP]{DP} J. Dziuba\'nski and M. Preisner, Remarks on spectral multiplier theorems on Hardy spaces associated with semigroups of operators,
\emph{Revista de la uni\'on Matem\'atica Argentina}, {\bf 50}
(2009), 201-215.

\bibitem[DVW]{DVW} B. Dahlberg, G. Verchota and T. Wolff, Unpublished manuscript.

\bibitem[He]{He} W. Hebisch, A multiplier theorem for Schr\"odinger
operators, {\it Colloq. Math.},  {\bf 60/61} (1990),  659-664.

\bibitem[H]{H} I. I. Hirschman, On multiplier transformations, Duke
Math. J., 26 (1959), 221-242.

\bibitem[LSV]{LSV} V. Liskevich, Z. Sobol and H. Vogt, On the $L^p$-theory of
$C_0$-semigroups associated with second-order elliptic operators. II,
J. Funct. Anal. 193 (2002), no. 1, 55�76.

\bibitem[MM]{MM} G. Mauceri and S. Meda, Vector-valued multipliers on stratified groups,
{\it Rev. Mat. Iberoamericana}, {\bf 6} (1990), 141-154.












\end{thebibliography}
\end{document}